\documentclass[epsfig,leqno,12pt]{amsart}

 \usepackage{eepic,epic}            

\theoremstyle{plain}

\theoremstyle{remark}

\def\AA{\mathbb{A}}

\def\SS{\mathbb{S}}
\def\TT{\mathbb{T}}
\def\Z{{\mathbb Z}}
\def\Q{{\mathbb Q}}
\def\C{{\mathbb C}}
\def\R{{\mathbb R}}

\def\al{\alpha}
\def\Aut{{\operatorname{Aut}}}

\def\gcd{{\operatorname{gcd}}}

\def\db[#1\db]{%
 \setbox0=\hbox{$#1$}\argwidth=\wd0
 \setbox0=\hbox{$\left[\box0\right]$}
    \advance\argwidth by -\wd0
 \left[\kern.3\argwidth\box0 \kern.3\argwidth\right]}

\begin{document}
\title[Growth partition functions\\ for \\cancellative infinite monoids 
]
{Growth partition functions\\ for \\ cancellative infinite monoids
}
\author{Kyoji Saito}
\address{ IPMU, university of Tokyo}
\thanks{}
 \centerline{} 
 \begin{abstract}
 We introduce the {\it growth partition function} 
 $Z_{\Gamma,G}(t)$ associated with any cancellative infinite monoid $\Gamma$ with a finite generator system
  $G$. 
 It is a power series in $t$ whose coefficients lie in integral Lie-like
  space $\mathcal{L}_{\Z}(\Gamma,G)$ in the configuration algebra associated with
  the Cayley graph $(\Gamma,G)$. 
  We determine them for homogeneous monoids admitting left greatest
  common divisor and right common multiple.
Then, for braid monoids and Artin
  monoids of finite type, using that formula, 
we explicitly determine their limit partition functions $\omega_{\Gamma,G}$. 
 \end{abstract} 

\maketitle   

\vspace{-0.7cm}
\tableofcontents

\vspace{-0.7cm}

\section{Introduction}\label{sec1} 

In a previous paper \cite[\S11.1.3]{S1}, 
we introduced the set
$\Omega(\Gamma,G)$ of {\it limit partition functions} (which
were called pre-partition functions there)
associated with a cancellative infinite monoid $\Gamma$ 
\footnote
{We call a semigroup with an identity element a {\it monoid}.  A
monoid is called {\it cancellative}  if the equality $aub\!=\!avb$ for
$a,b,u,v$ in the monoid implies $u\!=\!v$. If an equality $ab\!=\!c$ for
$a,b,c\!\in\! \Gamma$ holds in a cancellative monoid $\Gamma$, then $a$
(resp.\ $b$) is uniquely determined by $b,c$ (resp.\ $a,c$), which we
shall denote by $cb^{-1}$ (resp. $a^{-1}c$).
} 
with a fixed finite generator system $G$. 
In the present paper, using the same framework,
we introduce the {\it growth partition function}
$Z_{\Gamma,G}(t)$, which we have already studied without a
name (see (1.1) and following explanations in the next
paragraph, which relates the growth partition function with limit
partition functions). 
We determine the growth
partition functions for a class of homogeneous 
monoids which admit left greatest common divisors and right common multiples.
Then, using that result, we show that Artin
monoids of finite type, in particular braid monoids, up to possible
finite exceptions, are {\it simple accumulating} (i.e.\
$\Omega(\Gamma,G)\!=\!\{\omega_{\Gamma,G}\}$ for a single element
$\omega_{\Gamma,G}$), and then we  determine explicitly  the limit 
partition function $\omega_{\Gamma,G}$ for them  by solving an algebraic
equation arising from the denominator of their growth functions. 

\medskip
In the following, we briefly recall the definition of 
the set $\Omega(\Gamma,G)$ of limit partition functions 
associated with an infinite Cayley graph $(\Gamma,G)$
\footnote
{To be exact, we consider colored and oriented graph (see
Footnote 5). We shall sometimes assume further three
conditions {\bf H}, {\bf I} and {\bf S} on $(\Gamma,G)$ (see \S2), even
though they are unnecessary for the definitions of $\Omega(\Gamma,G)$
and $Z_{\Gamma,G}(t)$.
},
and recall in (1.1) the main formula in \cite{S1} for them. 
The main term of the formula is a proportion of two growth
functions, which we 
will call the {\it growth partition function} and denote by $Z_{\Gamma,G}(t)$.
  
An isomorphism class of a finite subgraph of $(\Gamma,G)$ is called  a
  {\it configuration}. The set of all configurations, denoted by
  $Conf(\Gamma,G)$, form a partial ordered semi-group by taking the disjoint
  union as the product. 
Consider the algebra $\AA[[Conf(\Gamma,G)]]:$ $\!=$the adic completion
  of the group ring $\AA\!\cdot\! Conf(\Gamma,G)$ with respect to the
  grading $\deg(S):\!=\!\#S$ for $S\!\in\! Conf(\Gamma,G)$, where $\AA$
  is a commutative coefficient ring. We can attach to it a topological
  Hopf algebra structure and call it the {\it configuration
  algebra}. It is also equipped with the classical topology if $\AA$ is $\R$ or $\C$. 
For any configuration $S\!\in\! Conf(\Gamma,G)$, let $\mathcal{A}(S)$ be
 the sum of  isomorphism classes of all subgraphs of $S$. Then,
 $\mathcal{M}(S)\!:=\!\log(\mathcal{A}(S))$ becomes a Lie-like element of the
  Hopf algebra, where we shall denote by $\mathcal{L}_{\AA}(\Gamma,G)$ the
  space of all Lie-like elements of  $\AA[[Conf(\Gamma,G)]]$.

Inspired by statistical mechanics, we call
  $\frac{\mathcal{M}(S)}{\#S}$ the {\it free energy} of $S$. 
We introduce the space $\Omega(\Gamma,G)$ of limit partition functions as
  the compact accumulation set (with respect to the classical topology)
 in $\mathcal{L}_\R(\Gamma,G)$ 
of the sequence
  of free energies $\frac{\mathcal{M}(\Gamma_n)}{\#\Gamma_n}$ for the
  balls  $\Gamma_n$ in $(\Gamma,G)$ of radius $n\!\in\!\Z_{\ge0}$
  centered at the unit $e$. Parallely, we introduce the {\it space
  $\Omega(P_{\Gamma,G})$ of  opposite series of the growth function}
$P_{\Gamma,G}(t)\!:=\!\sum_{n=0}^\infty\!\#(\Gamma_n)t^n$ 
(\cite[\S11.2.3]{S1}, see also \S2). Then, we obtain 
a natural surjective map:
\[
 \pi_\Omega\ :\ \Omega(\Gamma,G)\ \longrightarrow\ \Omega(P_{\Gamma,G})
\]
which is equivariant with certain actions $\tilde \tau_\Omega$ and $\tau_\Omega$ on $\Omega(\Gamma,G)$ and 
$\Omega(P_{\Gamma,G})$, respectively. Both actions are
transitive if $\Omega(\Gamma,G)$ is finite.
Therefore, the fiber of $\pi_\Omega$ over a point in 
$\Omega(P_{\Gamma,G})$ is an orbit of a finite
cyclic group $\Z/m_{\Gamma,G}\Z\!:=\!\mathrm{ker}(\langle\tilde\tau_\Omega\rangle\!\to\!\langle\tau_\Omega\rangle)$, called the {\it inertia group}.

The main formula of \cite[\S11.5 Theorem]{S1} states that 
 {\small
\begin{equation}
\vspace{-0.3cm}
\label{eq:trace}
 Trace^{[e]}\Omega(\Gamma,G)-E=\frac{m_{\Gamma,G}}{h_{\Gamma,G}}\sum_{x_i\in V(\Delta_{\Gamma,G}^{top})}
 A^{[e]}(x_i^{-1})\ \frac{P_{\Gamma,G}\mathcal{M}(t)}{P_{\Gamma,G}(t)}\Big|_{t=x_i}
\end{equation}
}
\!\!\!where $Trace^{[e]}\Omega(\Gamma,G)\!:=$the sum of limit partition function in a
fiber of $\pi_\Omega$ (=an orbit of the inertia group) over a point $[e]\!\in\!\Omega(P_{\Gamma,G})$, $E=$an error term
(conjecturally zero), $h_{\Gamma,G}=ord(\tau_{\Gamma,G})$,\footnote{See
\S2 for a definition of $\tau_{\Gamma,G}$.
The $h_{\Gamma,G}$ is called the {\it period}, characterized as the smallest integer s.t.\
$\Delta_{P_{\Gamma,G}}^{top}\!\big|(t^{h_{\small\Gamma,G}}\!\!-\!r_{\small\Gamma,G}^{h_{\Gamma,G}})$,
where $r_{\small\Gamma,G}$ is the radius of convergence of
$P_{\small\Gamma,G}(t)$ (see Footnote 4 {\bf A.}\ for a definition of $\Delta^{top}_{\Gamma,G}$). }
$V(\Delta_{\Gamma,G}^{top}):=$the set of zero loci of the
top-denominator polynomial
$\Delta_{\Gamma,G}^{top}(t)$ of $P_{\Gamma,G}(t)$ (see Footnote 4.{\bf A}),
$A^{[e]}(s)=$the numerator polynomial in $s$
of degree {\small$h_{\Gamma,G}\!-\!1$} of the opposite series 
indexed by $[e]\in\Omega(P_{\Gamma,G})$ 
\vspace{-0.2cm}
and, finally, 
\[\vspace{-0.2cm}
 P_{\Gamma,G}\mathcal{M}(t)\ :=\ \sum_{n=0}^\infty\mathcal{M}(\Gamma_n)t^n
\]
is a newly introduced growth function of Lie-like elements \cite[\S11.2.7]{S1}. 
Due to formula (1.1), 
we are interested in the ratio $\frac{P_{\Gamma,G}\mathcal{M}(t)}{P_{\Gamma,G}(t)}$, and  give it  a name: a {\it growth partition function} and denote it by  $Z_{\Gamma,G}(t)$. 

\medskip
In the following 1)-4), we contrast the growth partition function
$Z_{\Gamma,G}(t)$ with  limit partition functions in $\Omega(\Gamma,G)$.

\smallskip
1)   {\it Family over $\Omega(\Gamma,G)$ v.s.\  a single function with one variable $t$.}  

Limit partition functions  for $(\Gamma,G)$ are parametrized by a compact
set $\Omega(\Gamma,G)$, whereas  there is only one growth partition
function $Z_{\Gamma,G}(t)$ with one variable $t$, and data of
$Z_{\Gamma,G}(t)$ can be disclosed by specializing the variable $t$ to
special values $t_i$ at some zero loci of the {\it denominator}\! 
\footnote
{
Here, we are abusing the terminology  {\it denominator} of
$P_{\Gamma,G}(t)$ as follows. \\
{\bf A:}  In 1)-3), we consider the cases when the growth function
$P_{\Gamma,G}(t)$ belongs to $\C\{t\}_{r_{\Gamma,G}}$, where we put 
$\C\{t\}_r\!:=\!\{P(t)\!\in\!\C[[t]]\mid$ (i) $P(t)$ converges on the disc
$D(r)\!:=\!\{t\!\in\!\C\mid |t|\!<\!r\}$, and (ii) there exists, so
called, a denominator
polynomial $\Delta(t)$ in $t$ such that $\Delta(t)P(t)$ is holomorphic on a neighbourhood of $\overline{D}(r)$\}
for $r\in\R_{>0}$ (see \cite[S11.4 Def.]{S1}). 
Let
$\Delta_P(t)\!=\!\prod_i(t\!-\!x_i)^{d_i}$, where $x_i\!\in\!\C$ with
$|x_i|\!=\!r$ and $d_i\!\in\!\Z_{>0}$, be such denominator polynomial of minimal
degree. Then, by the top-denominator of $P(t)$, we mean
$\Delta_P^{top}(t)\!:=\!\prod_{i,d_i\!=\!d}\!(t\!-\!x_i)$ (where $d\!=\!\max\{\!d_i\!\}$).
If $\Omega(P)$ is finite, $\Delta_{P}^{top}(t)$ is a factor of
$t^{h}\!-\!r^{h}$ for
$h\!\in\!\Z_{>0}$ called the period \cite[\S11.3]{S1}.\\
{\bf B:} In 4), we consider the cases when the growth function is a
rational function  or a global meromorphic function in $t$ (they are
included in the above case {\bf A}).  Then the denominator of the growth
function means the denominator in usual sense, up to unit
factor. Obviously, $\Delta_P^{top}(t)$ 
is a factor of the denominator in this sense.
}
of the growth function $P_{\Gamma,G}(t)$. 
We do not know whether $Z_{\Gamma,G}(t)$ recovers the whole functions of $\Omega(\Gamma,G)$ or not. 
On the other hand, we shall see in the following 4) that
$Z_{\Gamma,G}(t)$ contains  ``new partition functions'' which may not be covered by the functions in $\Omega(\Gamma,G)$.

2)  {\it Completed coefficient field $\R$ v.s.\  small coefficient ring $\Z$.}  

We use the real number field $\R$ as the coefficient ring $\AA$ to
describe limit partition functions $\Omega(\Gamma,G)$, since they are defined by
classical limits of sequences of free energies whose
coefficients are in rational number field $\Q$, whereas the growth
partition function $Z_{\Gamma,G}(t)$ is defined as power series with
coefficients in integral lattice points
$\mathcal{L}_\Z(\Gamma,G)$ in the Lie-like space (see \S2 for the
lattice $\mathcal{L}_\Z(\Gamma,G)$). 

3) {\it Coefficients of partition functions.} 

For a reason in above 2), it is hard to determine explicit values of the coefficients
of limit partition functions in $\Omega(\Gamma,G)$ with respect to the
integral lattice basis. However, once
it is expressed using growth partition function, then the coefficients
appear explicitly by the substitution of the parameter  $t$ to some
special values: zero-loci  $x_i$ of the denominator polynomial of the growth
function, which are often ``calculable''.

4) {\it New partition functions.}  

In the above 1), 2) and 3),
$t$ is specialized at the zero-loci of denominators of the growth
function whose absolute values is the smallest (see {\footnotesize
Footnote 4.{\bf A}}). Let us call these partition functions tentatively
``old''. On the other hand, specialization of $Z_{\Gamma,G}(t)$ at other
zero-loci of the denominator (see {\footnotesize Footnote 4.{\bf B}})
give ``new partition functions'' in the sense that they satisfy the
kabi condition (see \cite[\S12, {\bf 2.\ Assertion}]{S1} and
Footnote 6.). The Galois group of the splitting field of the denominator polynomial acts on and mixes up old and new partition functions.

\medskip
The above 1), 2), 3) and 4) altogether seem to suggest that $Z_{\Gamma,G}(t)$ gives some structural insight on
partition functions, even though we do not yet understand the {\it
global phenomenon} described in 4) (see \cite[\S12, {\bf 2.}\ and {\bf
3.}]{S1} and \S4 Artin monoid of finite type). 

\medskip
Let us give an overview of the present paper.  

In \S2, we recall from \cite{S1} basic concepts and
notations on configuration algebras, introduce the space
$\Omega(\Gamma,G)$ of limit partition functions, and define the growth partition
function $Z_{\Gamma,G}(t)$.  We loosen a technical assumption in
\cite{S1} that $\Gamma$ is embeddable into a group to a weaker one, which
we call {\bf Assumption H}.
In \S3, we calculate the growth partition function for a class
$\mathcal{C}$ of cancellative homogeneous monoids which admit left greatest common divisors and
right common multiple. 
Finally in \S4, we show that an Artin monoid of finite type is {\it simple
accumulating}. Applying (1.1), 
we determine the unique limit partition function explicitly by a help of the
denominator polynomial of the growth function.

\section{Growth partition functions}

We recall basic notation and concepts (as minimal as possible) on configuration algebra on a Cayley
graph of a cancellative monoid (see \cite{S1} for details). Then, we
introduce the limit set $\Omega(\Gamma,G)$ of limit partition functions and
the growth partition function $Z_{\Gamma,G}(t)$.  


\smallskip
Let $(\Gamma, G)$ be the colored Cayley graph\footnote
{
Cayley graph $(\Gamma,G):=$ a graph whose vertex set is
$\Gamma$, and two vertices
$u,v\!\in\!\Gamma$ are connected by an edge if and only if $u^{-1}v$ or
$v^{-1}u\!\in\! G$. Each oriented edge $a-\!\!\!-\beta$ is labelled by the
element $\al^{-1}\beta$ in $G$ (called color and orientation). 
}
associated with a pair of  an infinite cancellative monoid $\Gamma$  and its finite generator system $G$. 
An isomorphism class denoted by $S\!=\![\SS]$ of a finite subgraph $\SS$ \footnote
{By a subgraph we mean a {\it full-subgraph}, i.e.\ two vertices are
connected by an edge in the subgraph if and only if they are connected in the Cayley
graph. Thus, an isomorphism of two subgraphs $\SS$ and $\TT$ is a bijection $\varphi$ of vertices such that, for $\al\!\in\!
G$ and $x,y\!\in\!\SS$, $x\al\!=\!y$ holds if and only if $\varphi(x)\al\!=\!\varphi(y)$ holds (see \cite[\S2.1]{S1})
}
 of
$(\Gamma,G)$ is called a {\it configuration}. The set of all
configurations (resp.\ connected configurations) is denoted by
$Conf(\Gamma,G)$ (resp.\ $Conf_0(\Gamma,G)$).  The set $Conf(\Gamma,G)$
has a monoid  structure generated by $Conf_0(\Gamma,G)$ by taking the
disjoint union as the product and the empty graph class $[\emptyset]$ as
the unit element. 
The completion 
$\mathbb{A}[[Conf(\Gamma,g)]]$ of the group ring $\mathbb{A}\cdot
Conf(\Gamma,G)$ 
with respect to the adic topology defined by the grading
$\deg(T)\!:=\!$ number of vertices of $T$  for $T\!\in\!
Conf(\Gamma,G)$ is called the {\it canfiguration algebra}, where
$\mathbb{A}$ is a commutative coefficient ring containing
$\mathbb{Q}$.
The configuration algebra is equipped with a topological Hopf algebra
structure, induced from the (higher) co-multiplications:
\vspace{-0.2cm}
{\small
\[
 \Phi_n:\  S\ \mapsto \ \sum_{S_1,\cdots,S_n\in Conf(\Gamma,G)}
\begin{pmatrix} S_1,\cdots,S_n\\ S\end{pmatrix}
S_1\otimes\cdots \otimes S_n
\vspace{-0.2cm}
\]
}
\noindent
for $n\!\in\!\Z_{\ge0}$ and $S\in Conf(\Gamma,G)$, where the coefficient is a combinatorial invariant, called the {\it covering
coefficient} (see \cite[\S2.4 \& \S4.1]{S1}). 

For $T\in Conf(\Gamma,G)$, let $\TT$ be a representative graph of $T$. Put
\[\begin{array}{lllll}
\!\!\! \mathcal{A}(T)\!&\!:=\!\!&\!\!\sum_{\SS\subset\TT}[\SS] = \text{\small the sum of 
  isomorphism classes of all subgraphs of $\TT$}.\!\!
\end{array}
\]
That is, $\mathcal{A}(T)\!=\!\sum_{S\in Conf}\!S A(S,T)$ for
$A(S,T)\!:=\!\#\AA(S,\TT)$ where $\AA(S,\TT)\!:=\!\{\SS\mid \SS\!\subset\! \TT\ \&\ [\SS]\!=\!S\}$.
Then $\mathcal{A}(T)$ is a {\it group-like element} 
in the Hopf algebra, i.e.\ $\Phi_n(\mathcal{A}(T))\!=\!\overset{n}\otimes
\mathcal{A}(T)$. In fact, this fact gives a
characterization of the Hopf algebra structure.
Thus, the logarithm $\mathcal{M}(T)\!:=\!\log(\mathcal{A}(T))$ for $T\! \in \!Conf(\Gamma,G)$  generate over
$\mathbb{A}$ a dense (w.r.t.\ the adic topology) submodule of the module
$\mathcal{L}_{\mathbb{A}}(\Gamma,G)$ of all Lie-like elements of $\AA[[Conf(\Gamma,G)]]$. However, they cannot form topological basis of 
$\mathcal{L}_{\mathbb{A}}(\Gamma,G)$, since $\mathcal{M}(T)\!=\!\#T\cdot [pt]\!+\!\cdots$ contains
low degree terms. Thus, we are lead to introduce a new (topological) $\AA$-basis
$\{\varphi(S)\}_{S\in Conf_0(\Gamma,G)}$ of $\mathcal{L}_{\mathcal{A}}(\Gamma,G)$
by the base change:
\begin{eqnarray}
\mathcal{M}(T) &=& \sum_{S\in Conf_0(\Gamma,g)}  \varphi(S)\cdot A(S,T),\\
\varphi(S)&=& \sum_{T\in Conf_0(\Gamma,g)}  \mathcal{M}(T)\cdot (-1)^{\#T-\#S}K(T,S),
\end{eqnarray}
\noindent
where $K(T,S)$
is a combinatorial constant $\in\! \Z_{\ge0}$, called kabi-coefficient, satisfying an inversion formula:$\sum_{U\in
Conf_0}(-1)^{\#U-\#S}K(S,U)A(U,T)=\delta(S,T)$ (\cite[\S7.3.1]{S1}). Further more, they satisfy
$A(S,T)\!=\!0$ if $S\!\not\le\! T$ and $K(T,S)\!=\!0$ if $T\!\not\le \!S$ or $\deg(S)\!\not\le\!\deg
(T)(\#G\!-\!1)\!+\!2$. 
In particular, $\varphi(S)$ consists only of terms of degree greater or
equal than $\deg(S)$ so that 
\vspace{-0.1cm}
\begin{equation}
 \mathcal{L}_{\mathbb{A}}(\Gamma,G) \ =\prod_{S\in Conf_0(\Gamma,G)}
\varphi(S)\cdot \mathbb{A}.
\vspace{-0.3cm}
\end{equation}
Regarding $\{\varphi(S)\}_{S\in Conf(\Gamma,G)}$ as  integral basis of
$\mathcal{L}_\AA(\Gamma,G)$, we put 

\noindent
$ \mathcal{L}_{\mathbb{B}}(\Gamma,G):=\prod_{S\in Conf_0(\Gamma,G)}
 \varphi(S)\cdot \mathbb{B}$
for any subalgebra $\mathbb{B}$ of $\mathbb{A}$.

Recall that for any $g\in \Gamma$, its length is defined by
\begin{equation}
\label{eq:length}
l(g):=\min\{n\in\Z_{\ge0}\mid \exists g_1,\cdots,g_n\in G \text{ s.t.\  } g=g_1\cdots g_n\}
\end{equation}
Note that $l(g)\!\ge\!\! d(g,e)\!\!:=\!$ the distance in the Cayley graph between\!
$g$ and the unit element $e$, but the equality may not hod in general. If $\Gamma$ is a group and $G\!=\!G^{-1}$,
the equality holds.

\medskip
\noindent
{\bf Definition.} We call a monoid $\Gamma$  {\it homogeneous with respect
to the generator system} $G$, if $l$ (2.5) is additive, i.e.\
$l(gh)\!=\!l(g)\!+\!l(h)$, or equivalently, if $\Gamma$ is presented by
homogeneous relations in $G$.

\medskip
Using $l(g)$, we define  a "ball" of radius $n\in\Z_{\ge0}$ centered at $e$ by
\begin{equation}
\label{eq:ball}
\Gamma_n:=\{g\in\Gamma\mid l(g)\le n\}.
\end{equation}
By an abuse of notation, we shall confuse the ball $\Gamma_n$ with its
isomorphism class in $Conf_0(\Gamma,G)$. 

\medskip
We recall definitions of
the spaces of limit partition functions and of opposite sequences,
and then state a Theorem on them. For the purpose, we specialize the
coefficient ring $\AA$ to the real number field $\R$. Here, we remark
that the configuration algebra $\R[[Conf(\Gamma,G)]]$ and the Lie-like
space $\mathcal{L}_\R(\Gamma,G)$ over $\R$ are also equipped with {\it classical topology}.

\medskip
\noindent
{\bf Definition.} {\bf 1.} The space of {\it limit partition function} of $(\Gamma,G)$ is
\[
 \Omega(\Gamma,G):=\begin{cases}
\text{
the accumulating set of free energies $\{\frac{\mathcal{M}(\Gamma_n)}{\#\Gamma_n}\}_{n\in\Z_{\ge0}}$} \\
\text{
in {\small $\R[[Conf(\Gamma,G)]]$} with respect to the
classical topology.
}
\end{cases}
\]
{\bf 2.} The space of opposite sequences for the growth sequence {\small $\{\#\Gamma_n\}$} is 
\[
 \Omega(P_{\Gamma,G}):=\begin{cases}
\text{
the accumulating set of polynomials {\small $\{\sum_{k=0}^n\!\!\frac{\#\Gamma_{n-k}}{\#\Gamma_n}s^k\}_{n\in\Z_{\ge0}}$}}\\
\text{
in {\small $\R[[s]]$} with respect to the
classical topology.
}
\end{cases}
\]
{\it Note.}  {\small Using formula (2.2), we see that the convergence of
$\frac{\mathcal{M}(\Gamma_n)}{\#\Gamma_n}$ on a subsequence of
$\{n\}_{n\in\Z_{\ge0}}$ is equivalent to the convergence of
$\frac{A(S,\Gamma_n)}{\#\Gamma_n}$ for all $S\in
Conf(\Gamma,G)$. Therefore, $\Omega(\Gamma,G)$ and $\Omega(P_{\Gamma,G})$ are closed subset of the Hilbert cubes 
$\prod_{S\in Conf_0(\Gamma,G)} \varphi(S)\!\cdot\! [0,1]$} and 
{\small $\prod_{n=0}^\infty s^n\!\cdot\! [0,1]$, respectively, 
so that} {\it $\Omega(\Gamma,G)$ and $\Omega(P_{\Gamma,G})$  are non-empty compact sets.}\footnote
{Further more, any element $\sum_{S\in
Conf_0}\varphi(S)\!\cdot\!a_S\!\in\! \Omega(\Gamma,G)$ satisfies a constraint $\sum_{S\in Conf_0}(-1)^{\#T-\#S}K(T,S)a_S\!=\!0$ for any
$T\!\in\! Conf_0$ (kabi condition \cite[\S11.1]{S1}). However, we shall not
discuss further on the condition in the present paper.}

\bigskip
\noindent
{\bf Theorem (\cite[11.2]{S1}).} {\it
Let $(\Gamma,G)$ be the Cayley graph of an infinite cancellative monoid
with a finite generator system $G$ satisfying {\bf Assumptions H}, {\bf
I'} and {\bf S} stated in the following proof.

1. The correspondence $\sum_S\varphi(S)\!\cdot\!a_S\mapsto \sum_k
s^k\!\cdot\!a_{\Gamma_k}$ defines a continuous surjective map:
\vspace{-0.2cm}
\[
 \pi_\Omega: \Omega(\Gamma,G) \rightarrow \Omega(P_{\Gamma,G})
\]

2. Define maps:  
\[
 \begin{array}{cccl}
 \tilde\tau_\Omega:&\mathcal{L}_\R\to \mathcal{L}_\R,&
 \sum_{S\in Conf_0}\!\varphi(S)\!\cdot\!a_S\mapsto
\frac{1}{a_{\Gamma_1}}\sum_{S\in Conf_0}\!\varphi(S)\!\cdot\!a_{S\Gamma_1}\!\!\!\!\!\!
\vspace{0.2cM}
\\
 \tau_\Omega:&\R[[s]]\to \R[[s]],&
\sum_{k=0}^\infty s^k\!\cdot\!a_k\mapsto
\frac{1}{a_{1}}\sum_{k=0}^\infty s^k\!\cdot\!a_{k+1},
\end{array}
\]
respectively, where i) their domains are restricted to the subspaces
$\{a_{\Gamma_1}\!\not=\!0\}$ and $\{a_1\!\not=\!0\}$, respectively, and ii)
$S\Gamma_1$ for $S\in Conf_0(\Gamma,G)$ means the isomorphism class of the graph
$\SS\Gamma_1\!=\!\cup_{\al\in\SS,\beta\in\Gamma_1}\al\beta$ for a
representative $\SS$ of the configuration $S$ (for the well-definedness of
$S\Gamma_1$, see {\bf 1.} in the following proof). 
Then, they induce continuous self-maps: 
\[
 \tilde\tau_\Omega:\Omega(\Gamma,G)\to\Omega(\Gamma,G)\quad \text{and}\quad
\tau_\Omega:\Omega(P_{\Gamma,G})\to\Omega(P_{\Gamma,G}), 
\]
respectively,
so that the  map $\pi_\Omega$ is equivariant with their
actions. That is, we obtain a commutative diagram:
\[\begin{array}{lccclll}
& \Omega(\Gamma,G)& \overset{\pi_\Omega}{\longrightarrow}& \Omega(P_{\Gamma,G})\\
& \tilde{\tau}_\Omega \downarrow\ \ &  & \tau_\Omega\downarrow\ \ \\
& \Omega(\Gamma,G)& \overset{\pi_\Omega}{\longrightarrow}&
  \Omega(P_{\Gamma,G})&.
\end{array}
\vspace{-0.3cm}
\]
}
\begin{proof} 
Theorem is already proven in \cite[\S11.2]{S1} under slightly stronger
 assumptions \cite[\S11.1 Assumption 1, \S11.2 Assumption 2.]{S1},
 which shall be replaced by {\bf H}, {\bf I'} and
{\bf S} given below. Therefore, in the following {\bf 1.} and {\bf 2}, we only explain
 new assumptions and sketch how they are used.  

\medskip
 {\bf 1.} In \cite[\S11.1 Assumption 1.]{S1}, we assumed
 that the monoid $\Gamma$ is embeddable in a group, say $\hat \Gamma$. 
Assumption 1.\ was used only to define the right action
 $\Gamma_1:Conf\to Conf,\ S\mapsto S\Gamma_1$. 
We shall replace Assumption 1.\ by the following weaker {\bf Assumption
H.}, which is sufficient to define the action of
 $\Gamma_1$, as we shall see in following {\bf Assertion A}. 
This weakening of the assumption shall be used when we study partition
 functions for a monoid of class $\mathcal{C}$ in \S3.

\medskip
\noindent
{\bf Assumption H.}  Assume the condition b) in next {\bf Assertion A}. holds. 

\medskip
We shall refer to this as the {\it homogeneity assumption on} $(\Gamma,G)$.

\bigskip
\noindent
{\bf Assertion A.}  {\it Let $\Gamma$ be a cancellative monoid
 generated by a finite set $G$. Then, in the following, a) implies b),
 and b) is equivalent to c).

a)  The monoid $\Gamma$ is embedded into a group, say $\hat\Gamma$.

b) {Let {\small $U_0,U_1,\cdots\!, U_n$ 
 and  $V_0,V_1,\cdots\!,V_n$
 ($n\!\in\!\Z_{\ge0}$)} be two sequences in $\Gamma$ such that every
 successive points {\small $U_{i-1},U_i$ and $V_{i-1},V_i$} for
 {\small $i\!=\!1,\!\cdots\!,n$} are connected by
 edges in {\small $(\Gamma,G)$}\! of the same label.  
If {\small $U_0\!=\!U_n$}  then {\small $V_0\!=\!V_n$}.
}

c) Any isomorphism $\varphi:\SS_1\!\simeq\!\SS_2$ between connected
 subgraphs of $(\Gamma,G)$
 induces an isomorphism $\hat \varphi:\SS_1\Gamma_1\!\simeq\!\SS_2\Gamma_1$
 such that  $\hat\varphi|_{\SS_1}=\varphi$.  
}
\begin{proof}
$a) \Rightarrow b)$:  
Regard $U_i$ and $V_i$ for $i=0,\cdots,n$ as elements in $\hat \Gamma$. Then
 $U_0\!=\!U_n$ implies that
 {\small
 $e=U_0^{-1}U_n\!=\!(U_0^{-1}U_1)(U_1^{-1}U_2)\cdots(U_{n-1}^{-1}U_n)$
 $\!=\!(V_0^{-1}V_1)(V_1^{-1}V_2)\cdots(V_{n-1}^{-1}V_n)\!\!=V_0^{-1}V_n$
} and $V_0\!=\!V_n$. 

\smallskip
$b) \Rightarrow c)$:  We need to show that i) a map $\hat\varphi:\SS\Gamma_1\to\SS_2\Gamma_1$ is
 well-defined by putting $\hat\varphi(\al\beta):=\varphi(\al)\beta$ for
 $\al\in\SS$ and $\beta\in G$, and ii) the map $\hat\varphi$ is an
 isomorphism of graphs. 

i) By definition of $\varphi$, if  $\al,\al\beta\in\SS_1$ and $\beta\in G$, then 
$\varphi(\al)\beta=\varphi(\al\beta)$. 

If
 $\al_1\beta_1=\al_2\beta_2\not\in \SS_1$ for $\al_1,\al_2\in \SS_1$ and
 $\beta_1,\beta_2\in G\cup G^{-1}$ where at least one of $\beta_1$ or
 $\beta_2$ belongs to $G$, then we need to show
 $\varphi(\al_1)\beta_1=\varphi(\al_2)\beta_2$.  Since $\SS_1$ is a
 connected graph, there is a sequence of points
 $U_1:=\al_1,U_2,\cdots,U_{n-1}:=\al_2$ which are successively adjacent
 to each other by an element of $G\cup G^{-1}$. Then, put $U_0:=U_1\beta_1$, $U_n:=\al_2\beta_2$ and $V_0:=\varphi(\al_1)\beta_1, V_1=\varphi(U_1),\cdots,V_{n-1}=\varphi(\al_2),V_n:=\varphi(\al_2)\beta_2$, we obtain two sequences satisfying the assumption in b). Then b) says that $V_0=V_n$ i.e.\ $\varphi(\al_1)\beta_1=\varphi(\al_2)\beta_2$.
Thus the map $\hat\varphi$ is well-defined. By applying the same argument for
 $\varphi^{-1}$, we see that $\hat\varphi$ is bijective. 

ii) It remains only to show that two distinct points $\al_1\beta_1$ and $\al_2\beta_2$ in $\SS_1\Gamma_1\setminus \SS_1$ is connected by an edge if and only if 
$\hat\varphi(\al_1\beta_1)$ and $\hat\varphi(\al_2\beta_2)$ are connected by the same labeled of edge.  This can be verified by comparing the sequence $\al_1\beta_1,\al_1,
\cdots, \al_2,\al_2\beta_2,\al_1\beta_2\gamma$ (here, $\al_1,\cdots,\al_2$ means a path in
 $\SS_1$ connecting the two points $\al_1$ and $\al_2$, and
 $\gamma:=(\al_2\beta_2)^{-1}\al_1\beta_1\in G$, by changing of the role
 of $\al_1,\beta_1$ and $\al_2,\beta_2$ if necessary) and the sequence
 $\hat\varphi(\al_1\beta_1),\hat\varphi(\al_1),\cdots,\hat\varphi(\al_2\beta_2),\hat\varphi(\al_2\beta_2)\gamma$. 
The condition b) says
 $\hat\varphi(\al_1\beta_1)=\hat\varphi(\al_2\beta_2)\gamma$, as
 desired.

 \smallskip
$c) \Rightarrow b)$: We show b) by induction on $n\!\in\!\Z_{\ge0}$, where the
 case $n\!=\!0$ is trivially true.
Let two sequences as in b) are given. If $U_i\!=\!U_j$ (resp.\ $V_i\!=\!V_j$)
 for $0\!\le\! i\!<\! j\! \le\! n$ and $(i,j)\!\not=\!(0,n)$, then
 by induction hypothesis, we have $V_i\!=\!V_j$ (resp. $U_i\!=\!U_j$). Then,
 applying the induction 
 hypothesis to the shorter sequences {\small
 $U_0,\!\cdots\!,U_i\!=\!U_j,\!\cdots\!,U_n$ and $V_0,\!\cdots\!,V_i=V_j,\!\cdots\!,V_n$},
 we obtain $V_0\!=\!V_n$. Thus, we may assume $U_0,\!\cdots\!,U_n$
 and $V_0,\!\cdots\!,V_n$ are mutually distinct except for $U_0\!=\!U_n$
 and  possible $V_0\!=\!V_n$.

Suppose we have $U_i\!=\!U_j\beta$ (resp.\ $V_i\!=\!V_j\beta$) for $\beta\!\in\! G$
 and $0\!\le\! i,j\!\le\! n$ such that $|i\!-\!j|\!\not=\!1$,
 $\{i,j\}\!\not=\!\{0,n\},\{0,n\!-\!1\}$ or $\{1,n\}$, then by applying the
 induction hypothesis to the sequences $U_i,\!\cdots\!,U_j,U_j\beta$ and
 $V_i,\!\cdots\!,V_j,V_j\beta$, we get $V_i\!=\!V_j \beta$
 (resp.\ $U_i\!=\!U_j\beta$). 

i)  Case 
$\beta\!:=\!U_{n-1}^{-1}U_n \!\in\! G$.  We have the natural isomorphism
 $\varphi:\SS_1\!=\!\{U_1,\cdots,U_{n-1}\}\simeq\SS_2\!=\!\{V_1,\cdots,V_{n-1}\},
 U_i\!\mapsto\! V_i\ (i\!=\!1,\!\cdots\!,n\!-\!1)$. The c) implies the
 existence of an isomorphism  $\hat\varphi\!:\!\SS_1\Gamma_1\!\simeq\!\SS_2\Gamma_1$. It then implies
 $\hat\varphi(U_0\!=\!U_n)\!=\!\hat\varphi(U_{n-1}\beta)\!=\!\varphi(U_{n-1})\beta\!=\!V_{n-1}\beta\!=\!V_n$. Since
 the vertex  $U_0$ is connected with $U_1$ by an edge, so is the
 vertex $\hat\varphi(U_0\!=\!U_n)\!=\!V_n$ with $\varphi(U_1)\!=\!V_1$. That is, in
 $\SS_2\Gamma_1$ two vertices $V_n$ and $V_0$ are connected with $V_1$
 by the same labeled edges, then the left cancellation implies $V_n\!=\!V_0$.

ii) Case $\beta:=U_n^{-1}U_{n-1}\in G$. Applying c) to the isomorphic
 graphs $\SS_1\!=\!\{U_0,\cdots,U_{n-2}\}$ and
 $\SS_2\!=\!\{V_0,\cdots,V_{n-2}\}$, isomorphism
 $\hat\varphi\!:\!\SS_1\Gamma_1\!\simeq\!\SS_2\Gamma_1$ implies
 $\hat\varphi(U_{n-1})\!=\!\varphi(U_{n-2})U_{n-2}^{-1}U_{n-1}\!=\!V_{n-1}$. On
 the other hand $U_0\!=\!U_n$ implies $U_{n-1}\!=\!U_0\beta$ and hence
 $\hat\varphi(U_{n-1})\!=\!\varphi(U_0)\beta\!=\!V_0\beta$. That is, two vetices
 $V_0$ and $V_n$ are connected with $V_{n-1}$ by edges of the same type
 $\beta$. Then the left cancellation by $\beta$ implies $V_0\!=\!V_n$.
\end{proof}
{\bf 2.} In \cite[\S11.2 Assumption 2.]{S1}, the following two were assumed.

\medskip
\noindent
{\bf Assumption I.}   Let $\mathbb{S}$ be a connected finite subgraph of
 $(\Gamma,G)$. If an equality  $\SS\Gamma_1=g\SS\Gamma_1$ for $g\in\hat
 \Gamma$ holds then $\SS=g\SS$ holds, where $\hat\Gamma$ is a group
 in which $\Gamma$ is embedded by Assumption 1.

\medskip
\noindent
{\bf Assumption S.} Define the set of dead elements\footnote
{
Since $\Gamma$ may not be a group and we do not assume $G=G^{-1}$, we should  note that our definition of dead elements is different from the followings
\begin{eqnarray*}
D_0(\Gamma,G)&:=&\{g\in\Gamma\mid l(h)\le l(g) \ \  \forall h\!\in\!\Gamma \text{ s.t.\  } h\!=\!g\al \text{ or } g\!=\!h\al \ \exists \al\!\in\! G\}.\\
D_1(\Gamma,G)&:=&\{g\in\Gamma\mid d(h)\le d(g) \ \  \forall h\!\in\!\Gamma \text{ s.t.\  } h\!=\!g\al \text{ or } g\!=\!h\al \ \exists \al\!\in\! G\}.
\end{eqnarray*}
}
by 
\begin{equation}
\label{eq:dead}
D(\Gamma,G):=\{g\in\Gamma\mid l(g\al)\le l(g) \ \  \forall \al\in G\}.
\end{equation}
 Then the ratio $ \frac{\#(\Gamma_n\cap D(\Gamma,G))}{\#\Gamma_n}$ tends to 0 as $n\to \infty$. 

\medskip
\noindent
{\it Note.} If $\Gamma$ is homogeneous with respect to $G$, then
 $D(\Gamma,G)\!=\!\emptyset$. Therefore, {\bf Assumption
 S.}\ is automatically satisfied.

\bigskip
Since we removed Assumption 1.\ that $\Gamma$ is embeddable into a
group $\hat \Gamma$, we need to reformulate {\bf I} in the following
 form {\bf I'} without using $\hat \Gamma$.

\bigskip
\noindent
{\bf Assumption I'.}   Let  $\SS$ and $\SS'$ be isomorphic finite
 connected  subgraphs of $(\Gamma,G)$. Then, any isomorphism
 $\hat\varphi\!:\!\SS\Gamma_1\!\simeq\! \SS'\Gamma_1$ induces $\hat\varphi|_{\SS}\!:\!\SS\simeq \SS'$.


\bigskip
Actually, Assumption {\bf I.} was used in \cite{S1} only
 once at  the proof of the following formula \eqref{eq:F}, 
which we will prove now by assuming only {\bf H}
 and {\bf I'} but not {\bf Assumption 1} and {\bf I}.

 \medskip
 \noindent
 {\it Formula. For $S\in Conf_0(\Gamma,G)$ and $n\in\Z_{>0}$,  we have }
\begin{equation}
\label{eq:F}
0\le A(S\Gamma_1,\Gamma_n)-A(S,\Gamma_{n-1})\le \#S\cdot \#(\dot \Gamma_n\cap D(\Gamma,G)).
\end{equation}
\noindent
{\it Proof of \eqref{eq:F}.}   The proof is parallel to  that of \cite[\S11.2.10]{S1}. For a sake of completeness of the present paper, we sketch it.

{\bf Assumption I'} implies that the  map $\cdot\Gamma_1:\AA(S,\Gamma_{n-1})\to
 \AA(S\Gamma_1,\Gamma_n)$ is injective. This implies the first
 inequality. 

On the other hand, any element
 $\mathbb{T}\in \AA(S\Gamma_1,\Gamma_n)$ is of the form $\SS\Gamma_1$
 for a graph $\SS\in \AA(S,\Gamma_n)$ ({\it Proof.}\ Fix $\SS_0$ with
 $[\SS_0]=S$. Put $\SS:=$ the image of $\SS_0$ by an  isomorphism
 $\SS_0\Gamma_1\simeq\mathbb{T}$. Then $\SS\subset\TT\subset\Gamma_n$).

If an element $\SS\Gamma_1\in \AA(S\Gamma_1,\Gamma_n)$ with $[\SS]=S$ is
 not in the $\Gamma_1$-image from $\AA(S,\Gamma_{n-1})$, i.e.\ $\SS\not\subset
 \Gamma_{n-1}$ then $\SS\cap\dot\Gamma_n\cap D(\Gamma,G)\not=\emptyset$.
 Let $\varphi:\SS_0\simeq\SS$ be an isomorphism. Choose points $d\in
 \SS\cap\dot\Gamma_n\cap D(\Gamma,G)$ and $s:=\varphi^{-1}(d)\in \SS_0$.
 Then, due to {\bf Assumption H.} and the connectedness of $\SS_0$, a pointed graph $(\SS,d)$ is uniquely determined (if it exists) as the isomorphic image of the pointed graph $(\SS_0,s)$, where the choice depends only on $(s,d)\in \SS_0\times (\dot\Gamma_n\cap D(\Gamma,G))$.  That is, the number of $\SS\Gamma_1\in \AA(S,\Gamma_n)$ with $\SS\not\subset \Gamma_{n-1}$ is at most $\#(S)\cdot \#(\dot\Gamma_n\cap D(\Gamma,G))$. 

This proves the the second inequality of \eqref{eq:F}. \qquad \qquad $\Box$

\medskip
Let us check how the inequality \eqref{eq:F} together with
 {\bf Assumption S.}\ imply the existence of the surjective map $\pi_\Omega$. First, we
 see easily
 that \eqref{eq:F} implies an inequality \cite[\S11.2.11]{S1}:
\[
 0\le A(\Gamma_k,\Gamma)-\#\Gamma_{n-k}\le
 \#(\Gamma_{k-1})\#(\Gamma_n\cap D(\Gamma,G))
\]
Then, one has $0\le \frac{A(\Gamma_k,\Gamma)}{\#\Gamma_n}-\frac{\#\Gamma_{n-k}}{\#\Gamma_n}\le
 \#\Gamma_{k-1}\frac{\#(\Gamma_n\cap D(\Gamma,G))}{\Gamma_n}
$, where {\bf Assumption S.}\ implies that the right hand side converges
 to 0 for any sub-sequence of  $\{n\}_{n\in\Z_\ge0}$ tending to infinity. Thus, the convergence of the
 first term to $a_{\Gamma_k}$ implies the convergence of the second term
 to $a_k$ such that $a_{\Gamma_k}=a_k$. This implies the map
 $\pi_\Omega$ is well defined. Surjectivity of the map $\pi_\Omega$ follows from
 the compactness of $\Omega(\Gamma,G)$: since for any sub-sequence
 $\{n\}_{n\in\Z_{\ge0}}$ tending to infinity such that
 $\frac{\#\Gamma_{n-k}}{\#\Gamma_n}$ converges for all $k\in\Z_{\ge0}$,
 we can choose sub-sequence of the subsequence such that
 $\frac{A(S,\Gamma_n)}{\#\Gamma_n}$ converges for all $S\in Conf(\Gamma,G)$.

In order to show  $\tilde\tau_\Omega(\Omega(\Gamma,G))\!\subset\!
 \Omega(\Gamma,G)$ and $\tau_\Omega(\Omega(P_{\Gamma,G)})\!\subset\!
 \Omega(P_{\Gamma,G)})$, 
using again the formula \eqref{eq:F} and {\bf  Assumption S.}, we show that
\[
\vspace{0.1cm}
\begin{array}{cll}
\tilde\tau_\Omega\big(\underset{m\to\infty}{\lim}\frac{\mathcal{M}(\Gamma_{n_m})}{\#\Gamma_{n_m}}\big)
 &=&
 \underset{m\to\infty}{\lim}\frac{\mathcal{M}(\Gamma_{n_m-1})}{\#\Gamma_{n_m-1}}\\
\\
 \tau_\Omega\big(\underset{m\to\infty}{\lim}\sum_{k=0}^{n_m}\frac{\#\Gamma_{n_m-k}}{\#\Gamma_{n_m}}s^k\big)
  &=&
  \underset{m\to\infty}{\lim}\sum_{k=0}^{n_m-1}\frac{\#\Gamma_{n_m-1-k}}{\#\Gamma_{n_m-1}}s^k .
\end{array}
\]
For the details of the proof, we refer to \cite[\S11.2]{S1}.

The equivariance of $\pi_\Omega$ with the actions $\tilde\tau_\Omega$
 and $\tau_\Omega$ is trivial since $\Gamma_k\Gamma_1=\Gamma_{k+1}$ for $k\in\Z_{\ge0}$.

This completes the proof of the Theorem.
\end{proof}
{\bf Question.} Do the conditions a), b) and c) in {\bf Assertion A.} equivalent?
Precisely, does b) imply a)? That is, do b) and c)
give characterizations of the embeddability of a monoid $\Gamma$ into a group?

\bigskip

Next in the remaining part of the present section, we introduce the {\it growth
partition functions} and discuss some of its descriptions. For the definition, we
do not need either {\bf Assumptions H.}, {\bf I'.}\ nor {\bf S}. Therefore,
until the end of this section 2, we assume only the
cancellativity on $\Gamma$.

Let us consider two growth series in the variable $t$: 
\begin{eqnarray}
P_{\Gamma,G}(t)&:=& \sum_{n=0}^\infty \#\Gamma_n \cdot t^n\\
P_{\Gamma,G}\mathcal{M}(t)&:=& \sum_{n=0}^\infty \mathcal{M}(\Gamma_n)\cdot t^n
\end{eqnarray}
where the first one is the usual growth function introduced by Milnor \cite{M} as an element of $\Z[[t]]$, and the second one is a growth series, introduced in \cite[(11.2.7)]{S1}  as an element in  $\mathcal{L}_{\Z[[t]]}(\Gamma,G)$.

\medskip
\noindent
{\bf Definition.} The {\it growth partition function} of
$(\Gamma,G)$ is the series
\begin{equation}
\label{eq:partition}
Z_{\Gamma,G}(t):=\frac{P_{\Gamma,G}\mathcal{M}(t)}{P_{\Gamma,G}(t)}
\end{equation}
Since the initial term $\#\Gamma_0$ of the growth function
$P_{\Gamma,G}(t)$ is $1$, the growth function is invertible in $\Z[[t]]$ so that  $Z_{\Gamma,G}(t)\in \mathcal{L}_{\Z[[t]]}(\Gamma,G)$.

\bigskip The following is an elementary remark.

\medskip
\noindent
{\it Assertion.} {\it The growth partition function has a development
\begin{equation}
Z_{\Gamma,G}(t)=\sum_{S\in Conf_0(\Gamma,G)}\varphi(S)\cdot Z_{\Gamma,G}(S,t).
\end{equation}
with respect to integral basis $\{\varphi(S)\}_{S\in Conf_0(\Gamma,G)}$, where
\begin{equation}
Z_{\Gamma,G}(S,t):= \frac{P_{\Gamma,G}A(S,t)}{P_{\Gamma,G}(t)}
\vspace{-0.4cm}
\end{equation}
\begin{equation}
P_{\Gamma,G}A(S,t):=\sum_{n=0}^\infty A(S,\Gamma_n) \cdot t^n
\end{equation}
(recall $A(S,\Gamma_n):=$the number of subgraphs in
$\Gamma_n$ isomorphic to $S$). 
}
\begin{proof}
Apply (2.2) for $T=\Gamma_n$, and, (2.10) together, we get
\begin{equation}
P_{\Gamma,G}\mathcal{M}(t) =\sum_{S\in Conf_0(\Gamma,G)}\varphi(S)\cdot P_{\Gamma,G}A(S,t).
\end{equation}
This together with (2.13) implies (2.12).
\end{proof}

It was shown \cite[10.6]{S1} that $P_{\Gamma,G}A(S,t)$ and $P_{\Gamma,G}(t)$ have same radius , say $r_{\Gamma,G}$, of convergence, so that the growth partition function converges at least in the radius $r_{\Gamma,G}$.  

\medskip
\noindent
{\bf  Conjecture 1.}  The growth partition function $Z_{\Gamma,G}(t)$
has  the radius of convergence larger than that $r_{\Gamma,G}$ of the
growth function $P_{\Gamma,G}(t)$.

\medskip
\noindent
{\bf  Conjecture 2.}  If the growth  function $P_{\Gamma,G}(t)$ is a
rational function in $t$, then the partition function coefficient
$Z_{\Gamma,G}(S,t)$ for any $S\in Conf_0(\Gamma,G)$  is a rational
function in $t$, whose order at infinity is bounded by
$L(S):=\min\{n\in\Z_{>0}\mid A(S,\Gamma_n)\not=0\}$.

\bigskip
\noindent
{\it Example.}
Let $F_f$ be a free group generated by
$G_f=\{g_1^{\pm1},\cdots,g_f^{\pm1}\}$ for $f\in\Z_{\ge0}$. 
The growth partition function for $(F_f,G_f)$ for $f\ge2$ is
\begin{eqnarray}
\label{eq:Free}
 \qquad
 Z_{F_f,G_f}(t)=\!\!\!\!\!
\sum_{\substack{S\in Conf_0(F_f,G_f), \\ d(S)\text{ even.}} } \!\!\!\!\!\!\!\!\! \varphi(S) t^{[d(S)/2]} 
+  \frac{2t}{1+t}   \!\!\!\! \sum_{\substack{ S\in Conf_0(F_f,G_f), \\ d(S) \text{ odd.}}} \!\!\!\!\!\!\!\!\!  \varphi(S) t^{[d(S)/2]} .\!\!\!\!\!\!\!\!\!\!\!\!\!\!\!\!\!\!\!\!\!\!\
\end{eqnarray}
where $d(S):=max\{d(x,y)\mid x,y\in S\}$ for  $S\in Conf_0(\Gamma,G)$. 

\begin{proof}
 For $n\ge[d(S)/2]$, the following formula holds \cite[\S11.1]{S1}:
\begin{eqnarray}
A(S,\Gamma_n) &=&\begin{cases}
  \frac{f(2f-1)^{n-[d(S)/2]}-1}{f-1} & \text{if $d(S)$ is even,}\\
  \frac{(2f-1)^{n-[d(S)/2]}-1}{f-1} & \text{if $d(S)$ is odd.}
\end{cases}
\end{eqnarray}
Thus, in view of the fact that $A(S,\Gamma_n)=0$ for $n<[d(S)/2]$, we calculate the growth function for $S$ as follows.
\begin{eqnarray}
\qquad \quad  P_{F_f,G_f}A(S,t)=
\begin{cases}
t^{[d(S)/2]}\frac{1+t}{(1-t)(1-(2f-1)t)} & \text{if $d(S)$ is even,}\\
t^{[d(S)/2]}\frac{2t}{(1-t)(1-(2f-1)t)} & \text{if $d(S)$ is odd.}
\end{cases}
\end{eqnarray}
In particular,
$\#\Gamma_n=\frac{f(2f-1)^n-1}{f-1}$ and $P_{F_f,G_f}(t)=\frac{1+t}{(1-t)(1-(2f-1)t)}$.

As a consequence, we obtain
\begin{eqnarray}
\qquad Z_{F_f,G_f}(S,t)=
\begin{cases}
t^{[d(S)/2]} & \text{if $d(S)$ is even,}\\
t^{[d(S)/2]}\frac{2t}{1+t} & \text{if $d(S)$ is odd.}
\end{cases}
\end{eqnarray}
Combining this with the formula (2.12), we obtain Formula \eqref{eq:Free}.
\end{proof}

\noindent
{\it Remark.} Specializing the growth  partition function at the two places  $t=1/(2f-1)=r_{F_f,G_f}$ and $t=1$ of poles of  $P_{F_f,G_f}(t)$, we obtain:
\begin{eqnarray}
\ \ \ \ \ \ \ \ \omega_{F_f,G_f}\!\!&\!=&\!\!\!\!\!\!\!\!\!\!
\sum_{\substack{S\in Conf_0(F_f,G_f), \\ d(S) \text{ even.}} }\!\!\!\!\!\!\frac{\varphi(S)}{(2f-1)^{[d(S)/2]} }
+  \frac{1}{f}\!\! \!\!\sum_{\substack{ S\in Conf_0 (F_f ,G_f), \\ d(S) \text{ odd.}}} \!\!\!\!\!\!\!	\frac{\varphi(S)}{(2f-1)^{[d(S)/2]}}\!\!\!\!\! \\
\ \ \ \ \ \ \ \ \  \omega_{F_f,1} \!\!&\!=&
\sum_{\substack{S\in Conf_0(F_f,G_f)} } \varphi(S).
\end{eqnarray}
where the first  formula coincides with the limit partition function for
$(F_f,G_f)$, which was already directly (without using $Z_{F_f,G_f}(t)$) calculated in \cite[\S11.1.9]{S1}. 


\section{Monoids of class $\mathcal{C}$}

We consider a class, which we call $\mathcal{C}$, of cancellative
homogeneous monoids with respect to finite generator system, admitting
conditions GCD$_l$ and CM$_r$.
We shall see that any configuration $S$ of a monoid of class $\mathcal{C}$ admits a
unique minimal representative, whose ``radius'' $L(S)$ is a numerical invariant of $S$.
Then, the growth partition function for the monoid of class
$\mathcal{C}$ is a sum of the main
term calculated by the invariant $L$ and the additional term coming
from dead elements.

\medskip
Let $\Gamma$ be a  cancellative monoid. Let us consider 
conditions on $\Gamma$:

\bigskip
\noindent
GCD$_l$:   For any two elements $u,v$ of $\Gamma$, there exists a unique
maximal common left divisor $\gcd_l(u,v)\in \Gamma$ of them. That is, 
$\gcd_l(u,v)|_lu$ and $\gcd_l(u,v)|_lv$, and if $w|_lu$ and $w|_lv$ for
$w\in\Gamma$ then $w|\gcd_l(u,v)$.

\smallskip
\noindent
CM$_r$:   For any two elements $u,v$ of $\Gamma$, there exists a common
right multiple of them. That is, there exists $w\in\Gamma$ such that
$u|_rw$ and $v|_rw$. In the other words, there exists $a,b\in \Gamma$ such that $au=bv$.

\bigskip
Note that uniqueness assumption in GCD$_l$,  in particular, asks that no
element in $\Gamma$ except for the unit element $e$ is invertible. In
particular, $\Gamma$  contains no non-trivial subgroups. 

\bigskip
\noindent
{\bf Definition.}  A  cancellative infinite homogeneous monoid
$\Gamma$ is called of class
$\mathcal{C}$ if it satisfies conditions GCD$_l$ and CM$_r$.

\bigskip
Let us state a general proporty of  monoids satisfying condition CM$_r$.

\bigskip
\noindent
{\bf Lemma 1.} {\it  Let $\Gamma$ be a cancellative monoid satisfying
condition CM$_r$. Then, for any finite generator system $G$ of $\Gamma$,
the Cayley graph $(\Gamma,G)$ satisfies {\bf Assumption H.}\
(recall \S2 for a definition). 
\begin{proof} 
Let {\small $U_0,U_1,\cdots\!, U_n$ 
 and  $V_0,V_1,\cdots\!,V_n$
 ($n\!\in\!\Z_{\ge0}$)} be two sequences in $\Gamma$ as in b) of \S2 Assertion A. Let us consider a common right multiple $W_0$ of $U_0$ and $V_0$. That
 is, there exists $A,B\in \Gamma$ such that $W_0=AU_0=BV_0$. We now
 compare two sequences {\small $AU_0,AU_1,\cdots\!, AU_n$ 
 and  $BV_0,BV_1,\cdots\!,BV_n$  ($n\!\in\!\Z_{\ge0}$)} in
 $\Gamma$. Let us show that the two sequences are the same. The initial
 terms have already the equality $AU_0=BV_0$. As induction hypothesis, assume
 $W_k:=AU_k=BV_k$ for a $k$ with $0\le k<n$. However, by the assumption
 on the sequences in b), two points $AU_{k+1}$
 and $BV_{k+1}$ are connected with $W_k$ by the same type edge. This implies that
 there exists $\al\in G$ such that either $AU_{k+1}=BV_{k+1}=W_k\al$ or 
$AU_{k+1}\al=BV_{k+1}\al=W_k$. In both cases, we get
 $AU_{k+1}=BV_{k+1}$. Thus we get finally $AU_{n}=BV_{n}$. 
If $U_0=U_n$, then  $BV_0=AU_0=AU_n=BV_n$. The left cancellation by $B$
implies $V_0=V_n$.
\end{proof}
\noindent
{\bf Corollary.} {\it Monoids of class $\mathcal{C}$ satisfies {\bf Assumption H}.}}

\bigskip
\noindent
{\bf Remark.}
If $\Gamma$ satisfies both CM$_r$ and CM$_l$ simultaneously, then, together with the
cancellativity, we know that $\Gamma$ is injectively embedded into its localization
group $\hat\Gamma$ 
of $\Gamma$ (\"Ore's criterion), implying {\bf Assumption H.}
%
We will observe in \cite{S4} that CM$_r$ alone together with cancellativity is sufficient
not only to get {\bf Assumption H.}, but leads to 
an embedding of $\Gamma$ into a ``homogeneous set'' $\hat \Gamma$ (which
may no longer have a group structure), where we define the growth partition
function for $\hat\Gamma$ (since for a
definition of limit partition functions and growth partition functions, group
structure is unnecessary \cite{S1}).

\bigskip
Let us return to the study  of monoids of class $\mathcal{C}$.

\bigskip
\noindent
{\bf Lemma 2.}  {\it Let $(\Gamma,G)$ be of class $\mathcal{C}$. Then, for any
configuration $S\!\in\! Conf_0(\Gamma,G)$, there exists a unique subgraph
$\SS_0$ of $(\Gamma,G)$ such that 
i) $[\SS_0]\!=\!S$ and ii) $\mathrm{gcd}_l(\SS_0)\!=\!e$. In particular, these imply
iii) $\mathrm{Aut}(\SS_0)\!=\!1$, and iv) for any subgraph $\SS$ with
$[\SS]\!=\!S$ we have $\SS\!=\!\gcd_l(\SS)\ \SS_0$.
}
\begin{proof}
Let $\SS$ and $\TT$ be any two presentative of the class $S$. That is, there
 is an isomorphism $\varphi:\SS\simeq\TT$. Choose any element
 $u\in\SS$. Let $w\in\Gamma$ be a common right multiple of $u$ and
 $\varphi(u)$. That is, there are $a,b\in \Gamma$ such that
 $w=au=b\varphi(u)$. Then, let us show that $a\SS=b\TT$ ({\it Proof.}
 It is sufficient to show that we have
 $ax=b\varphi(x)$ for any $x\in \SS$. But, this can be shown by 
 induction on the distance inside the graph $\SS$ of $x$ from $u$ by using
the cancellativity and  {\bf Assumption H.} of $\Gamma$.)

The uniqueness of the gcd$_l$ of elements of $a\SS\!=\!b\TT$ implies
 the equality: $a\ \gcd_l(\SS)\!=\!b\ \gcd_l(\TT)$. This implies the
 relation: $\gcd_l(\SS)^{-1}\SS\!=\!\gcd_l(\TT)^{-1}\TT$. This means that
 $\SS_0\!:=\!\gcd_l(\SS)^{-1}\SS$ does not depend on the choice of a
 representative $\SS$ of $S$. Thus, i), ii) and iv) are proven. 

Suppose that
 there is an automorphism $\varphi$ of $\SS_0$. Then, applying the same
 argument above, consider a common right multiple $au\!=\!b\varphi(u)$ for
 an element $u\!\in\! \SS_0$. Then the automorphism $\varphi$ is realized by
 $\SS_0\!\overset{\cdot a}{\sim}\! a\SS_0\!=\! b\SS_0\!\overset{/b}{\sim}\! \SS_0$. Again, the uniqueness
 of GCD$_l$ implies $a\!=\!b$ and $\varphi\!=\!1$.
\end{proof}
\noindent
{\bf Corollary.} {\it  Let $(\Gamma,G)$ be of type $\mathcal{C}$. Then, for any
configuration $S\in Conf_0(\Gamma,G)$, the automorphism group $\Aut(S)$
is trivial. In\! particular,\! $(\Gamma,G)$ satisfies {\bf Assumption I'}.\!
}

\bigskip
Let us call $\SS_0$ in Lemma 2.\ the {\it minimal
representative} of $S\!\in\! Conf(\Gamma,G)$. Using the minimal
representative, we introduce a numerical invariant for $S$, which is
used to present the growth partition function.

\medskip
\noindent
{\bf Notation.}  For  $S\in Conf_0(\Gamma,G)$, put  
\begin{equation}
\label{eq:L}
L(S):=\max\{ l(u)\mid u\in \mathbb{S}_0\}.
\end{equation}

\bigskip
\noindent  {\bf Formula.} The growth partition function for 
a pair of a monoid $\Gamma$ of class $\mathcal{C}$ and a finite generator
system $G$
is given by 
\begin{equation}
\label{eq:classc}
Z_{\Gamma,G}(t)=\sum_{S\in Conf_0(\Gamma,G)} \varphi(S) \ t^{L(S)}.
\end{equation}

\begin{proof}
Let us denote by $\mathbb{A}(S,\Gamma_n)$ the set of subgraphs of
 $\Gamma_n$ whose isomorphism class is equal to $S$ (recall \S2).  

\bigskip
\noindent
{\bf Lemma 3.} {\it  For $S\in Conf_0(\Gamma,G)$, let $\mathbb{S}_0$ be
 the minimal representative of $S$.
Then, for $n\in\Z_{\ge0}$, we  have a natural bijection: 
}
\begin{equation}
\label{eq:min}
\Gamma_n \simeq \mathbb{A}(S,\Gamma_{n+L(S)}),   \quad g\mapsto g\mathbb{S}_0.
\end{equation}
\begin{proof}
The correspondence is well-defined and is injective due to the uniqueness in Lemma 1.
Surjectivity is also clear from homogeneity of $(\Gamma,G)$, since if
 $\mathbb{S}\in \mathbb{A}(S,\Gamma_{n+L(S)})$ then, again by Lemma 1,
 we have $\mathbb{S}=\gcd_l(\mathbb{S})\cdot \mathbb{S}_0$, where
 $n+L(S)\ge \max\{l(u)\mid u\in \mathbb{S}\}=l(\gcd_l(\mathbb{S}))+L(S)$
 implies $n\ge l(\gcd_l(\mathbb{S}))$ and $\gcd_l(\SS)\in \Gamma_n$.
\end{proof}

\noindent
{\bf Corollary 1.}  {\it Under the same assumptions, we have the equality:
\begin{eqnarray}
P_{\Gamma,G}A(S,t) &=& t^{L(S)}\cdot  P_{\Gamma,G}(t) \\
Z_{\Gamma,G}(S,t)& =& t^{L(S)}.
\end{eqnarray}
\vspace{-0.3cm}
}
\begin{proof} We have $A(S,\Gamma_n)\!=\!0 \text{ if } n\!<\!L(S)$ and $ \#\Gamma_{n-L(S)}$ if $n\!\ge\! L(S)$.\!\!
\end{proof}
Corollary 1.\  together with (1.9) implies Formula \eqref{eq:classc}.
\end{proof}

As an application of Lemma 3., let us state about the map $\pi_\Omega$ (\S2).

\medskip
\noindent
{\bf Corollary 2.}  {\it Let $(\Gamma,G)$ be of class $\mathcal{C}$. Then, the map
$\pi_\Omega:\Omega(\Gamma,G)\to\Omega(P_{\Gamma,G})$ is a
bijection. 
}
\begin{proof} Generally, $\pi_\Omega$ is surjective under {\bf Assumption S.},
 which is automatically satisfied for homogeneous $(\Gamma,G)$ (recall
 {\it Note.}\ after \eqref{eq:dead}). 

For
 injectivity, we need to show that if the opposite polynomials  
 $X_n(P_{\Gamma,G}):=\sum_{k=0}^n\frac{\#\Gamma_{n-k}}{\#\Gamma_n}s^k$
  (recall \S2 Definition for $\Omega(\Gamma,G)$, \cite[\S11.2.2]{S1}) for a subsequence
 $\{n_m\}_{m\in\Z_{\ge0}}$ converges, then the free
 energies $\mathcal{M}(\Gamma_n)/\#\Gamma_n$ should converge also for
 the same subsequence $\{n_m\}_{m\in\Z_{\ge0}}$. Actually, using (2.1),
 for the convergence of the sequence $\mathcal{M}(\Gamma_n)/\#\Gamma_n$,
 it is sufficient to show the convergence of
 $A(S,\Gamma_n)/\#\Gamma_n$  for all $S\in Conf_0(\Gamma,G)$  for the same sequence $\{n_m\}_{m\in\Z_{\ge0}}$.

But it was shown (see
 \S2) that, under {\bf Assumptions I'.}\ and
 {\bf  S}, the sequence $A(\Gamma_k,\Gamma_n)/\#\Gamma_n$ and the sequence
 $\#\Gamma_{n-k}/\#\Gamma_n$ converge simultaneously for the same sequence $\{n_m\}_{m\in\Z_{\ge0}}$.
Thus,
 $A(S,\Gamma_n)/\#\Gamma_n=A(\Gamma_{L(S)},\Gamma_n)/\#\Gamma_n$ for
 $S\in Conf_0(\Gamma,G)$ \eqref{eq:min} for $n\ge L(S)$ converges for the same sequence $\{n_m\}_{m\in\Z_{\ge0}}$.
\end{proof}

Note that  the proof of Corollary 2. is independent from whether $\Omega(\Gamma,G)$ and/or  $\Omega(P_{\Gamma,G})$ is finite or not, and whether $P_{\Gamma,G}$ is a rational function or not.

\section{Artin monoids of finite type}

 Let $G$ be a finite set of letters and let
$M=(m_{\al,\beta})_{\al,\beta\in G}$ be a Coxeter matrix (i.e.\
$m_{\al,\al}\!=\!1$ for $\al\!\in\! G$ and $m_{\al,\beta}\!=\!m_{\beta,\al}\!\in\!
\Z_{\ge2}\!\cup\!\{\infty\}$ for $\al\!\not=\!\beta\!\in\! G$). Then, an Artin
monoid $\Gamma_M$ (\cite{B-S}) or a generalized braid monoid (\cite{D})
associated with the Coxeter matrix $M$ is a monoid defined by the positive homogeneous relations 
\[
\qquad \langle\al\beta\rangle^{m_{\al,\beta}}=\langle\beta\al\rangle^{m_{\al,\beta}}
\quad \text{for} \ \ \al,\beta\in G
\]
on the free monoid generated by the letters in $G$. Here, we denote by
$\langle\al\beta\rangle^m$  a word of alternating sequence of letters $\al$
and $\beta$ starting from $\al$ of length $m\!\in\!\Z_{\ge0}$. 
We shall refer to $G$ as the {\it standard generator system} for the Artin
monoid. 
An Artin monoid is called {\it of finite type}, if the associated Coxeter
group (i.e.\ the quotient group of $\Gamma_M$ divided by the relations
$\al^2\!=\!1$ for $\al\in G$) is finite. Indecomposable Artin monoids of finite type
are classified into types $A_l$
($l\ge1$), $B_l$ ($l\ge2$), $D_l$ $(l\ge4)$, $E_6$, $E_7$, $E_8$, $F_4$,
$G_2$, $H_3$, $H_4$ and $I_2(p)$ ($p\ge3$). 
A {\it braid monoid} $B(n)^+$ of n-strings is isomorphic
to an Artin monoid of type A$_{n-1}$. 

\medskip
\noindent
{\it Assertion. An Artin monoid of finite type belongs to the class $\mathcal{C}$. }
\begin{proof} 
Clearly, any Artin monoid is homogeneous by the definition.

 It is shown (\cite{B-S},\cite{D}) that an Artin monoid is an cancellative infinite monoid, satisfying conditions GCD$_l$ and
GCD$_r$. If, further,  it is of finite type, it they satisfies LCM$_l$ and
LCM$_r$, and hence CM$_l$ and CM$_r$.  
\end{proof}

\noindent
{\it Note.} An Artin monoid is known to be embeddable into its group
(Paris \cite{P}) so that it satisfies {\bf Assumption H}. However, we shall not
use this result in the present paper. See also {\bf Remark 3.}\ at the end of the paper.

\medskip
 Recall (\cite{S2},\cite{S3}) that the  growth function for an Artin monoid $\Gamma_M$ with
 respect to the standard generator system $G$ is given by 
\begin{equation}
P_{M}(t) :=P_{\Gamma,G}(t)=  \frac{1}{N_M(t)},
\end{equation}
where 
\begin{equation}
\label{eq:denom}
N_M(t):= \sum_{J\subset G}(-1)^{\#J} t^{\deg(\Delta_J)}.
\end{equation} 
Here the summation index $J$ runs over all subsets of $G$ such that the
restriction $M|_J\!:=\!(m_{ij})_{i,j\in J}$ is a Coxeter matirx of
finite type, and $\Delta_J$ is the fundamental element in the monoid
 $\Gamma_{M|_J}$ so that $\deg(\Delta_J)\!=$length of the longest
element in the associated Coxeter group (\cite{B-S}).

\smallskip
For an indecomposable Artin monoid of finite type $\Gamma_M$, the following (1), (2) and (3) are
conjectured \cite{S2}.

(1)  $\tilde N_M(t):=N_M(t)/(1-t)$ is an irreducible polynomial over
$\Z$, 

(2)  there are $\#G-1$ distinct real roots on the interval $(0,1)$, and 

(3) the smallest real root on the
 interval (0,1), say $r_{\Gamma,G}$,  of $N_M(t)\!=\!0$ is strictly
 smaller than the absolute value of any other root. 

\medskip
 Actually, conjectures are affirmatively solved for types $A_l,
 B_l\!=\!C_l$  and $D_l$ for $l\!\le\!30$ and $E_6,E_7,
 E_8,F_4,G_2,H_3,H_4$ and $I_2(p)$ ($p\!\ge\!3$) by a help of
 computer. Conjecture (3)   is affirmatively solved by
 Kobayashi-Tsuchioka-Yasuda \cite{K-T-Y} for the 
types $A_l$, $B_l$ and $D_l$ for $l\!\ge\! M$ for some $M$, where the author
 still expect that $M=1$.

\medskip
Conjecture 3 implies 
$\Delta_{P_{\Gamma,G}}^{top}(t)=t-r_{\Gamma,G}$, where $r_{\Gamma,G}$ is the
radius of convergence of the series $P_{\Gamma,G}(t)$.
As a consequence, the period $h_{\Gamma,G}$ is equal to 1 and
$\#\Omega(P_{\Gamma,G})\!=\!1$ (recall Footnote 3)  except for
possible finite exceptions in types $A_l, B_l$ or $D_l$ with
$30\!<\!l\!<\!M$.  Together with \S3 Corollary 2. to Lemma 3, this implies
$\#\Omega(\Gamma,G)=1$, that is, $(\Gamma,G)$ is {\it simple accumulating}
in the terminology of \cite[\S11.1]{S1}. Let us denote by $\omega_{\Gamma,G}$
the single element of $\Omega(\Gamma,G)$. Then, since we have $\omega_{\Gamma,G}=Trace^{[e]}\Omega(\Gamma,G)$ where the class $[e]$ denotes
the single element in $\Omega(P_{\Gamma,G})$, this limit element is now calculable
from the growth partition function \eqref{eq:classc}
by the use of a formula \eqref{eq:trace}. Let us describe the other
terms in \eqref{eq:trace}. 

$E=$ a sum of terms depending the root of
$\delta:=(t^{h_{\Gamma,G}}-r_{\Gamma,G}^{h_{\Gamma,g}})/\Delta^{top}_{\Gamma,G}$

\ \ \ $= 0$ \quad since $\deg(\delta)=0$  (see \cite[(11.3.6), (11.5.6)]{S1}). 

$h_{\Gamma,G}=\#(\Omega(P_{\Gamma,G}))=1$ (Conjecture 3 in \cite{S1}, solved by \cite{K-T-Y}). 

$m_{\Gamma,G}=$ covering
sheet number of $\Omega(\Gamma,G)\to\Omega(P_{\Gamma,G})$ 

\qquad \ $=1 $ \ \ (see \S3 Lemma 3, Corollary 2).

$A^{[e]}(s)\!=\!$ a polynomial in $s$ of degree $h_{\Gamma,G}\!-\!1$ with
a constant term 1 

\qquad\ \  \ $=1 $, \ \ since $h_{\Gamma,G}-1=0$ (see  \cite[\S11.3.2]{S1}).

\medskip
Finally, substituting in $Z_{\Gamma_M,G}(t)$ (3.23) the single summation index $r_{\Gamma_M,G}$: the
smallest real root of the equation $N_M(t)=0$, we arrive at the goal
formula of the
present paper.

\bigskip
\noindent
{\bf Theorem.}
{\it Artin monoid of finite type, except for finite possible exceptions
in types $A_l, B_l$ or $D_l$, is simple accumulating. That is,
$\#\Omega(\Gamma_M,G)\!=\!1$.  The limit partition function is given by}
\begin{equation}
\omega_{\Gamma_M,G} = Z_{\Gamma_M,G}(r_{\Gamma_M,G})= 
\sum_{S\in Conf_0(\Gamma_M,G)} \varphi(S) \ r_{\Gamma_M,G}^{L(S)},
\end{equation}
{\it where $r_{\Gamma_M,G}$ is the smallest real root in the interval (0,1)
of the denominator polynomial} \eqref{eq:denom}.

\bigskip
\noindent 
{\bf Remark 1.}
Let $M$ be an indecomposable Coxeter matrix of finite type. Let us consider the set $\Delta:=\{\al\in\C\mid \tilde N_M(\al)=0\}$ of roots of $N_M(t)=0$, and, for $\al\in\Delta$,  put
\begin{equation}
\label{eq:}
\omega_{M,\al}:=Z_{\Gamma_M,G}(t)|_{t=\al}=\sum_{S\in Conf_0(\Gamma_M,G)} \varphi(S) \ \al^{L(S)}.
\end{equation}
It was shown (\cite[\S11.4, {\bf 4.} Assertion.]{S1}) that each $\omega_{M,\al}$ belongs to the Lie-like space $\mathcal{L}_{\C,\infty}$ at infinity.
Then assuming Conjecture (1) in \cite{S2},\cite{S3}, the Galois group of
the splitting field of $N_M(t)$ acts transitively on the set $\Delta$,
inducing also a transitive action on the set
$\{\omega_{M,\al}\}_{\al\in\Delta}$ of limit partition functions. In particular, the action mixes up
the limit partition function $\omega_{\Gamma,G}$ with the other functions.
We do not know the meaning of this action.

\medskip
{\bf 2.}  Above Theorem is valid not only for Artin monoids of finite
type but for any monoid $(\Gamma_,G)$ of type $\mathcal{C}$ whose growth
function belongs to $\C\{t\}_{r_{\Gamma,G}}$ and  has period
$h_{\Gamma,G}$ equal to 1 (i.e. $r_{\Gamma,G}$ is the unique pole of
$P_{\Gamma,G}$ on the circle $|t|=r_{\Gamma,G}$). 

\medskip
{\bf 3.}
Artin monoids of non-finite type do not belong to the class
$\mathcal{C}$, since they do not satisfy CM$_r$. 
However,  we conjectured in \cite[\S3 Conjecture 3]{S3} that Artin
monoids of affine type have the period $h_{\Gamma,G}$ equal to 1. 
Thus, it seems to be interesting to ask whether \S3 {\bf Lemma 2.}\ of
the present paper holds for Artin
monoids (in general or, in particular, of affine type) or not.

\medskip

{\bf Acknowledgment.} This work was supported by World Premier
International Reseach Center Initiative (WPI Initiative), MEXT, Japan.

\end{document}